\documentclass[amslatex, 12pt]{amsart}

\input amssym.def
\input amssym.tex

\input xy
\xyoption{all}

\def\codim{\mathop{\rm codim}}
\def\depth{\mathop{\rm depth}}
\def\dim{\mathop{\rm dim}}
\def\projdim{\mathop{\rm projdim}}
\def\reg{\mathop{\rm reg}}
\def\rank{\mathop{\rm rank}}

\def\Su{{\mathcal S}}

\def\C{{\mathcal C}}

\def\iP{{\mathbb P}}

\def\im{{\mathfrak m}}

\def\o{{\omega}}

\newtheorem{thm}{Theorem}[section]
\newtheorem{lemma}[thm]{Lemma}
\newtheorem{cor}[thm]{Corollary}
\newtheorem{pro}[thm]{Proposition}
\newtheorem{example}[thm]{Example}
\newtheorem{remark}[thm]{Remark}
\newtheorem{defn}[thm]{Definition}
\newtheorem{notn}[thm]{Notation}
\newtheorem{blank}[thm]{}

\newcommand{\bt}{\begin{thm}}
\newcommand{\et}{\end{thm}}
\newcommand{\blem}{\begin{lemma}}
\newcommand{\elem}{\end{lemma}}
\newcommand{\bco}{\begin{cor}}
\newcommand{\eco}{\end{cor}}
\newcommand{\bp}{\begin{pro}}
\newcommand{\ep}{\end{pro}}
\newcommand{\bex}{\begin{example}}
\newcommand{\eex}{\end{example}}
\newcommand{\brm}{\begin{remark}}
\newcommand{\erm}{\end{remark}}
\newcommand{\bd}{\begin{defn}}
\newcommand{\ed}{\end{defn}}
\newcommand{\bnot}{\begin{notn}}
\newcommand{\enot}{\end{notn}}
\newcommand{\bblank}{\begin{blank}}
\newcommand{\eblank}{\end{blank}}

\newcommand{\bib}{\bibitem}

\newcommand{\sms}{\setminus}
\newcommand{\beqn}{\begin{eqnarray*}}
\newcommand{\eeqn}{\end{eqnarray*}}
\newcommand{\beq}{\begin{eqnarray}}
\newcommand{\eeq}{\end{eqnarray}}
\newcommand{\been}{\begin{enumerate}}
\newcommand{\eeen}{\end{enumerate}}
\newcommand{\lrar}{\longrightarrow}

\newcommand{\Q}{\mathbb Q}
\newcommand{\N}{\mathbb N}
\newcommand{\Z}{\mathbb Z}
\newcommand{\iS}{{\overline S}}

\usepackage{mathrsfs}
\usepackage{graphicx}
\usepackage{chancery}
\usepackage{bbding}
\usepackage{pifont}
\usepackage{fancybox}
\usepackage{palatino}
\usepackage{color}
\usepackage{picinpar}

\newcommand{\mptm}{\fontfamily{ptm}\fontsize{10}{12}\selectfont}

\begin{document}

\title[Canonical module]{  On the canonical module of  toric surfaces in $\iP^4$}

\author{Clare D'Cruz}
\address{Chennai Mathematical Institute, Plot H1, SIPCOT IT Park
Padur PO, Siruseri 603103, India}
\email{clare@cmi.ac.in}




\maketitle

\mptm
\section{Introduction}

Let $S$ be a semigroup ring and let $k[S]$ denote the semigroup ring over a field $k$. We will assume that $k[S]$ is equidimensional. 
In this paper we are interested in the Cohen-Macaulayness of the canonical module $\o_{k[S]}$. If $k[S]$ is Cohen-Macaulay then there is nothing to prove. Hence we will assume that $k[S]$ is not Cohen-Macaulay. 

Interest in the canonical module 
 was motivated by a paper of Mumford \cite{mumford} where he
showed that  for an ample line bundle ${\mathcal L}$ on a normal
surface $\Su$, $H^{1}(\Su, {\mathcal   O}_{\Su}\otimes {\mathcal
  L}^{-1})=0$.  He remarks that 
this is false if $\Su$ doesn't satisfy $S_{2}$ ({\it i.e.} is not
Cohen-Macaulay) and asks if  the $S_{2}$ condition is
sufficient. The first counter-example was given in \cite{ara}.

The study of the canonical module of toric
ideals has been of recent interest. 
In this paper we show that the canonical module of a toric surface
$\Su$ in $\iP^4$ is Cohen-Macaulay, in other words the $S_2$-ification
of $\Su$ always arithmetically Cohen-Macaulay. Hence if $\Su$ satisfies $S_2$, then
$H^{1}(\Su,{\mathcal   O}_{\Su}\otimes {\mathcal L}^{-1})=0$. 
As a consequence  the canonical module $\o_{k[S]}$ is Cohen-Macaulay. 

Of recent is also the question of the  Castelnuovo-Mumford regularity. The famous conjecture of Eisenbud and Goto states that for a non-degenerate irreducible closed scheme over an algebraically closed field $\reg({\mathcal{S}}) 
\leq \deg({\mathcal{S}}) - \codim({\mathcal{S}})$. 
The Eisenbud-Goto conjecture for the bound on the regularity of toric varieties of codimension two was proved in 
  \cite{strum-peeva}. Theorem~\ref{thm-main}(\ref{thm-main-reg}) that for a toric surface in $\iP^4$, any 
  surface linked to a toric surface by a complete intersection ideal, (i.e., it is in the CI-liaison class of a toric surface in $\iP^4$) satisfies the Eisenbud-Goto conjecture.

It is well known that the  Castelnuovo-Mumford  canonical module of a
reduced curve is always two. More generally, if the dimension of the
ring is at least two then  the depth of the canonical
module is at least two. Hence it is natural to ask if the canonical
module of a surface satisfies $S_3$, i.e., it is Cohen-Macaulay. 

We first show that if $S$ is a standard semi-group ring then $\o_{k[s]}$ is always Cohen-Macaulay.  This does not hold true in general, even for a toric surface in any  projective space. For toric surface we prove:

\bt
\label{thm-main}
Let $S$ define a toric surface of in $\iP^4$. Let $I$ be the defining ideal of 
$k[S]$ and    $J$ an ideal in the CI-liaison class of $I$. Then the following hold true
\been

\item
 $R/ I \cong \oplus_{n \in \Z} H^0(\Su, {\mathcal O}_{\Su}(n))$. 

\item
 $R/ J  \cong \oplus_{n \in \Z} H^0( (R/J)^{\sim}, {\mathcal O}_{(R/J)^{\sim}}(n))$. 
 
 \item
$\o_{R/I}$ and $\o_{R/J}$ is Cohen-Macaulay.

\item
The $S_2$-ification of $R/I$ and $R/J$ are Cohen-Macaulay.

\item
\label{thm-main-reg}
$reg(R/J) \leq deg(R/J) - \codim(R/J) + 1$.

\item
$\reg(\o_{R/I})= \reg(\o_{R/J})=3$.
 \eeen
\et

We also answer a  question of Goto and Watanabe in
\cite{goto-w}. Let $S$ be an affine semi-group ring. If $S$ is $S_2$
then is  $k[S]$ is Cohen-Macaulay? This was shown to be false in 
\cite{trung-hoa}. 
In this paper we show that if Goto and Watanabe's conjecture holds true if $k[S]$ defines a
toric surface in 
$P^4$. 

\noindent
{\bf Acknowledgement:}  This work began while the author was visiting Institut de Math{\'e}matiques
CNRS \& Universit{\'e} Paris 6. The author thanks Supported by Minist{\`e}re de la recherche for financial support and Universit{\'e} Paris 6 for local hospitality.  The author also wishes to  Marc Chardin who motivated this question and for his encouragement during the preparation of this paper. 

\section{Preliminaries}
The following notations in this section are from \cite{trung-hoa}.

\bnot
Let  $S$ be  an affine semi-group 
ring in ${\N}^n$,  $n \in \N $. 
\been
\item
 $G_S$ is  the additive group 
generated by $S$ in ${\Z}^n$.

\item
$\iS = \{ x \in G : n g \in S \mbox{ for some } n \geq 0
\}.$

\item
For  $1 \leq i \leq n$,
$S_{(i)} = \{ x \in G_S: x + y \in S \mbox{ for all } y \in S \cap F_i\}$ and $S_i = S - S_{(i)}$.

\item
$C(S) = \mbox{convex rational polyhedral cone spanned 
by $S$ in ${\Q}^n$}$.

\item
$m = \#\{ \mbox{$r-1$-dimensional faces of $C(S)$ } \}$. 

\item
$\pi(S) := $ the simplicial complex of nonempty subsets  $J$ of $[1,m]$ such that 
$\cap_{i \in J} S \cap F_i \not = (0) $. 

\item
For any subset  $J\subseteq [1,m]$, let 
$G_J := \cap_{i \not \in J} S_i  \sms \cup_{j \in J} S_j$

\item
For any subset $J \subseteq [1,m]$, let $\pi_J$ be the simplicial complex of nonempty subsets 
$I$ of $J$ such that $\cap_{i \in I} S \cap F_i \not = (0) $. 

\item
$S' = \cap_{i=1}^n S_i$.
\eeen
\enot

\bd
We say that an affine semi-group is standard if
\been
\item
$\iS = G_S \cap  {\N}^n$, 

\item
Let $1 \leq i,j \leq n$. $S_{(i)} \not = S_{(j)}$ for $i \not = j$.

\item
$\rank_{\Z}G_{S_{(i)}} = r-1$.

\eeen
\ed

\bblank {\bf Hochster Transform} \cite{hochster}
We describe the Hochster transform:
Consider  the vector space $V$ generated by 
$S \cup \{0\} \in {\Q}^n$.  

There exists linear functionals $L_1, \ldots, L_m$ such that 
$L_i(x) \geq 0$ for all $ x \in V$. 
Then $L := {[L_1, \cdots, L_m]}^{t} : {\Q}^n \lrar {\Q}^m$. Then $L(V)
\subseteq {\Q}^m$. By multiplying $L_i$ by a suitable positive
multiple we can assume that $L(S) \subseteq {\N}^m$. 
\eblank

\blem
Let $S$ be an affine semi-group ring. Then
\been
\item
$L(S) \cong S$;

\item
${\overline{L(S)}} = G_{L(S)} \cap {\N}^m$

\item
$L(S)$ is a standard affine semigroup ring.
\eeen
\elem

\bd
Let  $S$ be a standard affine semi-group in ${\N}^r$. 
\been
\item
Let $\pi(S)$ be the simplicial complex of non-empty subsets $J$ of
$[1,m]$  with the property that $\cap_{i \in J} S \cap F_i \not =
(0)$.

\item
Let $\pi_J$ be the set simplicial complex of non-empty subsets $I$ of $J$ with
the property $\cap_{i \in I} S \cap F_i \not = (0)$.
\eeen  
\ed

Let $k[S]$ be  a homogeneous semigroup ring which is homomorphic image of a polynomial ring 
$R = k[x_0, x_1, \ldots, x_n]$. 

\bd
The canonical module of $S$ is 
\beqn
{\o}_S := Ext^c(S, \o_R)
\eeqn
where $c = \codim~S = dim~S - dim~R$. 
\ed


\section{Main Result}
Our main result is the Cohen-Macaulayness of the canonical module of the 
projective toric variety. Since $k[S']$ is the $S_2$-ification of $S$ and $\omega_s \cong \omega_{S'}$, it is enough to show that $k[S']$ is Cohen-Macaulay. We state two results  of Trung and Hoa which is describes the cohomology groups of $k[S']$.

\bt
\label{trhoa}
\cite[Corollary~3.7]{trung-hoa}
For all $i < r$,
\beqn
      H^i(k[S']) 
\cong \bigoplus_{J \not \in \pi(S); \#J \leq m-2}
       k[G_J] \otimes_k {\tilde H}_{i-2}(\pi_J;k)
\eeqn
\et

\blem
\label{lemma1}
\cite[Corollary~3.4]{trung-hoa}
Let $S$ be a  semi-group in ${\N}^{r}$ of rank $r$.
Then
\beqn
      H^1(k[S']) =0.
\eeqn
\elem

In view of Lemma~\ref{lemma1}, it is clear that  the $k[S']$ is Cohen-Macaulay if and only if 
 $H^i(\o_{k[S']}) = 0$ for $\dim~k[S] \geq i \geq 2$. 
Hence if $S$ defines a toric surface in $\iP^4$, then $k[S']$ is Cohen-Macaulay if and only  $H^2(\o_{k[S']}) = 0$.

\blem
\label{lemma2}
Let $S$ be a  standard semi-group in ${\N}^{r}$ of rank $r$. Then $k[S']$ is Cohen-Macaulay. 
\elem
\proof 
Applying Lemma~\ref{lemma1}, we only need to show that 
 $H^i(k[S']) =0$ for all $2 \leq i \leq r-1$.
After applying the Hochster Transform we can assume that 
\beqn
         \pi(S) 
=        \{\phi, 
\cup_{i=1}^r \{i\},  
\cup_{1 \leq i < j \leq r} \{ij\} \cdots \cup
\cup_{i=1}^{r}  \{1 \ldots {\hat{i}} \cdots n\} \}.
\eeqn
Hence
\beqn
\{J \not \in \pi(S); \#J \leq r-2\} = \phi.
\eeqn
Now applying Theorem~\ref{trhoa} we get  $H^i(k[S'])=0$ 
for $1=2, \ldots,r-1$. 
\qed

\blem
\label{theorem2}
Let $S$ be an affine semigroup in ${\N}^3$. Assume that $m=3$. 
Then $k[S']$ is Cohen-Macaulay. 
\elem
\proof If $m=3$, then by applying the Hochster transform we can assume that $S$ is a standard semigroup ring of rank $3$ in $\N^3$. Hence by Lemma~\ref{lemma2},
$k[S']$ is Cohen-Macaulay.
\qed

\blem
\label{theorem3}
Let $S$ be an affine semigroup in ${\N}^3$. Assume that $m=4$. 
Then $k[S']$ is Cohen-Macaulay. 
\elem
\proof Applying the Hochster transform we can assume that our semigroup has generators are in 
  ${\N}^4$ and are of the form
\beqn
\begin{array}{lllllr}
s_1 & = &  [ 0          & s_{12}   & s_{13}   & 0      ]  \\
s_2 & = &  [ 0          & 0            & s_{23}   & s_{24} ]  \\
s_3 & = &  [ s_{31} & 0            & 0            & s_{34} ]  \\
s_4 & = &  [ s_{41} & s_{42}   & 0            & 0      ]  \\
s_5 & = &  [ s_{51} & s_{52}   & s_{53}   & s_{54} ]  \\
 \end{array} 
\eeqn
where $s_5$ is in the convex hull of the cone generated by $s_1$, $s_2$, $s_3$ and $s_4$. Hence we can express 
\beqn
\begin{array}{llll}
s_5 =  q_1 s_1 +   q_2 s_2 +  q_3 s_3 + q_4 s_4
\end{array}
\eeqn
where $q_i \in Q$ and  $q_i \geq 0$ with 
$q_1 + q_2 + q_3 + q_4 = 1$.
Without loss of generality we can assume that any $s_1$, $s_2$ and $s_3$ are lineraly dependent.

We have   
\beqn
\pi(S) &=& \{ \phi ,
\{1\},\{2\},\{3\},\{4\},\{12\},\{14\},\{23\},\{34\} \}\\
J &=& \{13\}, \{24\}.
\eeqn
We apply Theorem~\ref{trhoa}. It is enough   show that $G_{J} = \emptyset$ 
for $J =  \{13\}$ as the proof for $J =\{24\}$ is similar. 

Now 
\beqn
G_{13} = (S_2 \cap S_4) \backslash (S_1 \cup S_3).
\eeqn  
 
Let $x \in S_2 \cap S_4$. Then we can express $x$ in two different ways:
\begin{alignat}{2}
     x  
&=  \label{one}n_1 s_1 +  (n_2 - n_2') s_2 +  (n_3 - n_3') s_3 + n_4 s_4 +  n_5 s_5  & \mbox{ as } x \in S_2\\
&= \label{two}
 (m_1- m_1') s_1 +  m_2 s_2 +  m_3 s_3 +  (m_4 - m_4') s_4 +  m_5 s_5 
 & \mbox{ as } x \in S_4.
\end{alignat}
If  $(n_2 - n_2') \geq 0$, then $x \in S_3$. If 
  $(n_3 - n_3') \geq 0$, then $x \in S_1$. 
If both are negative then since $x \in S_4$ we consider (\ref{two}).
Once again, if $m_1 -m_1' \geq 0$, then $x \in S_3$, otherwise if $m_4- m_4' \geq 0$, then $x \in S_1$. 
Therefore if  atleast one of the four integers
$n_2- n_2'$,  $n_3 - n_3'$, $m_1 - m_1'$ and $m_4 - m_4'$ is non-negative, we are done. 

{\bf Claim 1:}  At least one of the four integers
$n_2- n_2'$,  $n_3 - n_3'$, $m_1 - m_1'$ and $m_4 - m_4'$  is non-negative.

Suppose  not, i.e., all of them are negative. 
We express $x$ component-wise. 
From (\ref{one}) and (\ref{two}) we get
\beq
\label{three}
     (n_3 - n_3') s_{31} + n_4 s_{41}          + n_5 s_{51}
&=& m_3 s_{31}           + (m_4 - m_4') s_{41} + m_5 s_{51}\\
    \label{four}
    n_1 s_{12}           + n_4 s_{42}          + n_5 s_{52} 
&=& (m_1 - m_1') s_{12}  + (m_4  - m_4')s_{42} + m_5 s_{52}\\\
     \label{five}
     n_1 s_{13}          + (n_2 - n_2') s_{23} + n_5 s_{53} 
&=& (m_1  - m_1')s_{13}  + n_2  s_{23}         + m_5 s_{53}\\
    \label{six}
   (n_2 - n_2') s_{24}   + (n_3 - n_3') s_{34} + n_5 s_{54}
  &=& m_2 s_{24}         + m_3 s_{34}          + m_5 s_{54}.
  \eeq
From (\ref{four}) and (\ref{six}) we have:
\begin{alignat}{2}
\label{eight}
        (m_1 - m_1'-n_1) s_{12}  + (m_4  - m_4'-n_4)s_{42} 
      &=(n_5- m_5) s_{52}; 
 &\qquad     & \\
\label{nine}
     (n_2 - n_2'-m_2) s_{24}   + (n_3 - n_3'-m_3) s_{34} 
  &=    (m_5  - n_5)s_{54}; 
 &\qquad & .
\end{alignat}
Since all the four integers $m_1 - m_1'-n_1$, $m_4  - m_4'-n_4$, 
$n_2 - n_2'-m_2$ and $n_3 - n_3'-m_3$ are negative, we conclude from (\ref{eight})that 
 $(n_5- m_5) s_{52} \leq 0$ and from (\ref{nine}) that $(m_5  - n_5)s_{54} \leq 0$. 

\noindent 
{\bf Claim 2:} $(n_5 - m_5)s_{52} \not = 0$: 
If $(n_5 - m_5) s_{52} = 0$, then since $m_1 - m_1'-n_1$ and  $m_4  - m_4'-n_4$ are negative, we have $s_{12} = s_{42} = 0$ and consequently $s_{52} = q_1 s_{12} + q_4 s_{42} = 0$. This implies that $m \leq 3$ which leads to a contradiction. 

\noindent 
{\bf Claim 3:} $(m_5  - n_5)s_{54} \not =  0$. 
If $(m_5 - n_5) s_{54} = 0$, then since $n_2 - n_2'-m_2$ and  $n_3  - n_3'-m_3$ are negative, we have $s_{24} = s_{34} = 0$ and consequently $s_{54} = q_2 s_{24} + q_3 s_{34} = 0$. This implies that $m \leq 3$ which leads to a contradiction. 

\noindent
Claim~2 implies that $n_5  < m_5$ and Claim~3 implies that $m_5 < n_5$ which is not possible. This proves Claim~1.

Hence  $S_2 \cap S_4 \subseteq S_1 \cup S_3$, and $G_{13} = (S_2 \cap S_4) \backslash (S_1 \cup S_3) = \emptyset$. This implies  that $k[G_{13}]=0$. Similarly
we can show that $k[G_{24}]=0$. 
\qed

\blem
\label{theorem4}
Let $S$ be a standard affine semigroup in ${\N}^3$. Let $m=5$. 
Then   $k[S']$ is Cohen-Macaulay.
\elem
\proof Applying  the Hochster transform we  can assume that our semigroup is in   $\Z^5$ and has generators
\beqn
\begin{array}{lllllllrl}
s_1 & = &  [ 0       & 0        & s_{13}   & s_{14}   & s_{15}  &] & \\
s_2 & = &  [ s_{21}  & 0        & 0        & s_{24}   & s_{25}  &] & \\
s_3 & = &  [ s_{31}  & s_{32}   & 0        & 0        & s_{35}  &] & \\
s_4 & = &  [ s_{41}  & s_{42}   & s_{43}   & 0        & 0       &] & \\
s_5 & = &  [ 0       & s_{52}   & s_{53}   & s_{54}   & 0       &] &.\\
 \end{array} 
\eeqn
We may  assume that $s_1$, $s_2$, $s_3$ are linearly independent and that 
\beqn
s_4 &=& a_1 s_1 + a_2 s_2 + a_3 s_3 \\
s_5 &=& b_1 s_1 + b_2 s_2 + b_3 s_3 
\eeqn
where $a_i, b_i \in \Q$.
Now
\beqn
\pi(S) 
&=& \{ J \subseteq [1,5] : \cap_{i \in J} S \cap F_i \not =(0) \} \\
                &=& \{ \emptyset,  \{1 \},  \{ 2 \}, \{3 \},   \{ 4 \} , \{5 \},  \{12   \}, \{15 \}, \{23 \}, \{ 34 \}, \{45 \}  \}\\
 \{ J \not \in \pi(S) \}
&=& \{ \{13\},  \{14\},  \{24\}, \{25\},  \{35\}, \{123\}, \{124\}, 
       \{125\}, \{134\}, \{135\}, \{145\}, \{234\},\{235\}, \{345\} \}.
\eeqn

We will show that $G_{13} = \emptyset $. The proof is similar for all 
$J \in \{ J |  \in  \pi(S) \mbox{ and } \#J = 2 \}$.

If $\#J = 3$, for example let $J = {135}$, then 
\beqn
G_{135} = (S_2 \cap S_4) \backslash (S_1 \cup S_3 \cup S_5).
\eeqn
Since $x \in S_2 \cap S_4$, following the proof of $G_{13} = \emptyset$, we can show that $x \in S_1 \cup S_3 \subseteq S_1 \cup S_3 \cup S_5$. 
 Similarly one can show that $G_J = \emptyset$  for all 
 $J \in\{ J |  \in  \pi(S) \mbox{ and } \#J = 3 \}$.

Consider $G_{13}$. 
\beqn
G_{13} = (S_2 \cap S_4 \cap S_5) \backslash (S_1 \cup S_3).
\eeqn
Let $x \in S_2 \cap S_4 \cap S_5$, then $x \in S_2 \cap S_4$. 
Hence we can write
\mptm
\begin{alignat}{2}
      x  
&= \label{oone} (n_1 - n_1') s_1 +  (n_2  - n_2') s_2 +  n_3 s_3 + n_4 s_4 +  n_5 s_5  & \mbox{ as } x \in S_2\\
      &=  \label{ttwo}m_1 s_1 +  m_2 s_2 +  (m_3 - m_3') s_3 +  (m_4 - m_4') s_4 +  m_5 s_5 & \mbox{ as } x \in S_4.
 \end{alignat}
Writing $x$ componentwise, from (\ref{oone}) and (\ref{ttwo}) 
\begin{alignat}{2}
     (n_2 - n_2') s_{21} + n_3 s_{31}          + n_4 s_{41}\label{ffour}
&=  m_2 s_{21}          + (m_3 - m_3') s_{31} + (m_4 - m_4') s_{41}
&&
\\
     n_3 s_{32}          + n_4 s_{42}          + n_5 s_{52}\label{ffive}
&= (m_3 - m_3') s_{32}  + (m_4 - m_4') s_{42} + m_5 s_{52}
&&
\\
    (n_1 - n_1') s_{13}  + n_4 s_{43}          + n_5 s_{53}\label{ssix}
&= m_1 s_{13}           + (m_4 - m_4') s_{43} + m_5 s_{53}
&&
\\
    (n_1 - n_1') s_{14}  + (n_2 - n_2') s_{24} + n_5 s_{54}\label{sseven}
&= m_1 s_{14}           + m_2 s_{24}          + m_5 s_{54}
&&
\\
    (n_1 -n_1') s_{15}    + (n_2 - n_2') s_{25} + n_3 s_{35}\label{eeight}
&= m_1 s_{15}           + m_2 s_{25}          + (m_3 - m_3') s_{35}
&&
\end{alignat}
{\bf Claim 1:} All the for integers $n_1 - n_1'$, $n_2 - n_2'$, $m_3 - m_3'$, $m_4 - m_4'$ cannot be negative. 

If all of them are negative,
then from (\ref{ffive}) and (\ref{sseven}) we get
\begin{alignat}{2}
\label{seventeen}
   (m_3 - m_3' - n_3) s_{32} + (m_4 - m_4' - n_4)s_{42} 
&= (n_5 - m_5) s_{52}; \\
\label{eighteen}
    (n_1 -n_1' - m_1) s_{14}    + (n_2 - n_2'- m_2) s_{24} 
&=  (m_5 - n_5)  s_{54}. 
 \end{alignat}
By our assumption  all  the integers $m_3 - m_3' - n_3$, $m_4 - m_4' -n_4$, $n_1 - n_1' - m_1$ and  $n_2 - n_2' - m_2$ are all negative. 

\noindent
{\bf Claim 2:} $s_{52} \not = 0$.: 
If $s_{52} = 0$, then since  $m_3 - m_3' - n_3 <0$ and $m_4 - m_4' -n_4<0$, we have $s_{32}= s_{42}=0$.  
This implies that $m=4$.

\noindent
{\bf Claim 3:} $s_{54} \not = 0$.
If $ s_{54} = 0$, then from (6),  $s_{14} = s_{24} = 0$.   This implies that 
$m=4$.
 
 Claim~2  implies that $n_5 - m_5 < 0$ while 
Claim~3  implies that $m_5 - n_5 < 0$. This leads to a contradiction and proves Claim 1. 

If $n_1 - n_1' \geq 0$, then from (\ref{oone}) we have $x \in S_3$.
If $n_2 - n_2' \geq 0$, then once again from (\ref{oone}) we have $x \in S_1$.
If $m_4 - m_4' \geq 0$, then 
we have $x \in S_3$.

\noindent
{\bf Claim 4:} The following is not possible:
$n_1 < n_1'$, $n_2 < n_2'$, $m_3 \geq m_3'$ and $m_4 < m_4'$.

Suppose not. From (\ref{eeight}) we have:
\begin{alignat}{1}
 (n_1 -n_1' - m_1) s_{15}    + (n_2 - n_2'-m_2) s_{25} 
&= (m_3 - m_3'-n_3) s_{35} .
\end{alignat}
Since $n_1 -n_1' - m_1 <0$ and $n_2 - n_2'-m_2<0$, we have 
$m_3 - m_3'-n_3 \leq  0$. 

{\bf Claim 5} $ s_{35} \not = 0$. 

If not, then,  then $s_{15} = s_{25}=0$. This implies  that  $m = 4$. 
This proves Claim 5. 

Hence $m_3 - m_3'-n_3 <  0$. From (\ref{seventeen}) we have
\begin{alignat}{1}
 (m_3 - m_3' - n_3) s_{32} +(m_4 - m_4' - n_4)s_{42} 
 &= (n_5 - m_5) s_{52} 
\end{alignat}
where $m_3 - m_3' - n_3<0$ $m_4 - m_4' -n_4<0$.
Using the arguement in Claim 5 we have  $n_5 < m_5$. 
From (\ref{eighteen}) we have $m_5 < n_5$. 
 This leads to a contradiction.
This proves Claim 6.
\qed

\section{Applications}
We would like to list some interesting applications.

\bt
\label{thm-final}
Let $S$ be either a standard semigroup  or let $S$ define a toric surface $\Su$ in ${\iP}^4$. 
Then 
\been
\item
$k[S']$ is Cohen-Macaulay;

\item
$\o_S$ is Cohen-Macaulay.
\eeen
\et
\proof For a standard semi-group ring (1)  has been proved in Lemma~\ref{lemma2}. For a toric surface in $\iP^4$, (1) follows from 
Lemma~\ref{lemma1}, Lemma~\ref{lemma2}, 
Lemma~\ref{theorem3} and Lemma~\ref{theorem4}.

Since $\o_{k[S]}  = \o_{k[S']}$ and $\o_{k[S']}$ is Cohen-Macaulay, 
$\o_{k[S]}$ is Cohen-Macaulay. 
\qed

\bt
\label{reg-omega}
Let $S$ be a semi-group ring which defines a  toric surface $\Su$ in ${\iP}^4$. Then 
the Castelnuovo-Mumford regularity of $\o_S$ is $3$.
\et
\proof Since $S$ is a domain,  $S$ is reduced and consequently  $S'$ is reduced.  Therefore 
$H^0(\Su', {\mathcal O}_{\Su'}(-n) = 0$ for $n >0$ and  
$H^0(\Su', {\mathcal O}_{\Su'}) \not = 0$. Also since $S'$ is Cohen-Macaulay,
 $H^1_{\im}(\Su', {\mathcal O}_{\Su'}) = 0$. Now 
 $H^3_{\im}(\o_S)_n = H^0(\Su',  {\mathcal O}_{\Su'}(-n)$. Hence 
 $\reg(\o_S) = 0+ 3 = 3$. 
\qed

\blem
\label{depth}
Let $V$ be a toric variety in ${\iP}^n$ of dimension at least two and codimension 
two. Then the depth of its coordinate ring $ S = R/I_{V}$ is at least  $n-2$. 
\elem
\proof This follows from \cite[Theorem~2.3]{strum-peeva}, that the 
$\projdim~S \leq 2^{\codim(I)} -1 = 2^2 - 3$.  Hence by the Auslander Buchsbaum formula, $\depth~S \geq (n+1) -3 = n-2$. 
\qed

\blem
Let $V$ be a toric variety in $\iP^n$ $(n \geq 4)$ of codimension 
two. Then 
$H^i(S_{\im}) = 0$ for 
$i=0, 1$ and 
hence $S \cong \oplus_{\mu}H^{0}(\Su ,{\mathcal O}_{\Su}(\mu ))$. 
\elem
\proof Follows from Lemma~\ref{depth}.
\qed


%

%

%

%

\noindent
{\bf Proof of Theorem~\ref{thm-main}:} By Lemma~\ref{depth}, $H_{\im}^0(R/I) = H_{\im}^1(R/I) =0$. Hence (1) follows. 

It is enough to prove (2), (3), (5) and (6) for an ideal $J$ linked to $I$ by a complete intersection ideal as it will follow for any ideal $J$ in the CI-liaison class. 

Now, let  $I$ and   $J$ be linked by a complete intesection ideal $I_{{\mathcal C}}$. Then we have the short exact sequence 
\beq
\label{sect-three-one}
0 \lrar \o_{R/I} \lrar R/ I_{\mathcal{C}} \lrar R/J \lrar 0
\eeq 
Since  $H^i(R/I_{\mathcal{C}}) = 0$ for $i \leq 2$, we have for $i=0,1$
\beqn
H^i(R/ J) \cong H^{i+1}( \o_{R/I})= 0
\eeqn
as $\o_{R/I}$ is Cohen-Macaulay. 
This proves (2). 

Interchanging $I$ and $J$ in (\ref{sect-three-one}), we have the long exact sequence of cohomology modules and for $i=1,2$
\beqn
H^i(R/ I) \cong H^{i+1}( \o_{R/J})= 0 \hspace{.5in} \mbox{[from  (1)]}.
\eeqn
Since $H^1(R/I_{\mathcal{C}}) = 0$, $H^1(\o_{R/J}) = 0$.  Hence $H^{i}( \o_{R/J})= 0$ for all $i \leq 2$. This proves (3). 
 
The Cohen-Macaulayness of the $S_2$-ification of $R/I$ has been proved in Section~3. Since $\o_{R/J}$ is Cohen-Macaulay and the $S_2$-ification of $R/J$ is $\o_{\o_{R/J}}$ which is Cohen-Macaulay. This proves (4).

\beqn
\label{sect-three-one-long}
0 \lrar H^2_{\im}(R/J) 
\lrar H^3_{\im}(R/ \o_{R/I})
\lrar H^3_{\im}(R/I_{\mathcal C}) 
\lrar H^3_{\im}(R/J) 
\lrar 0 
\eeqn
Since $R/ I_{\mathcal{C}}$ is a complete intersection ideal and  $\o_{R/I}$ is Cohen-Macaulay, we have  $H_{\im}^0(R/I) = H_{\im}^1(R/I) =0$. Hence (2) follows. 

We have already shown that $\o_{R/I}$ is Cohen-Macaulay. To show that $\o_{R/J}$ is Cohen-Macaulay interchange the role of $I$ and $I$ in the exact sequence (\ref{sect-three-one}) and use (2). 
The $S_2$-ification of $R/J$ is  $\o_{\o_{R/J}}$ which is Cohen-Macaulay since $\o_{R/J}$. Hence (4) follows. 

It remains to prove (5). By \cite[Theorem~7.3]{strum-peeva}
\beqn
reg(R/I) &\leq& deg(R/I) - \codim (R/I) + 1\\
         & < & deg(R/ I_{\C}) - \codim (R/I) + 1\\
         &=&   d_1 + d_2 - 2 + 1 \\
         &=&   d_1 + d_2 -1.
\eeqn
From the exact sequence
\beqn
0 \lrar \o_{R/J} (5 -  d_1 - d_2)\lrar R/I_{\C} \lrar R/I \lrar 0
\eeqn
This gives rise to the long exact sequence 
\beqn
0 \lrar H^2_{\im}(R/I)
  \lrar H^3_{\im} (\o_{R/J})(5 - d_1 - d_2)
  \lrar H^3_{\im} (R/ I_{\C})
  \lrar  H^3_{\im}(R/I)
  \lrar 0.
\eeqn
We get 
\beqn
        reg(\o_{R/J}) 
&\leq&  d_1 + d_2 -2  - d_1 - d_2 + 5 =    3. 
\eeqn
\qed

We now give a positive answer to a question of Goto and Watanabe for toric surfaces in $\iP^4$. 

\bt
Let $S$ be a semigroup ring which defines a toric surface in $\iP^4$. Then 
$S = S'$ if and only if $k[S]$ is Cohen-Macaulay.
\et
\proof
If $S=S'$, then $k[S] = k[S']$ and $k[S']$ is Cohen-Macaulay by
Theorem~\ref{thm-final}. Conversly if $k[S]$ is Cohen-Maculay, then $k[S]_{\im}$ is Cohen-Macaulay for every maximal ideal $\im$ of $R$ and 
$End(\o_{k[S]_{\im}}) \cong k[S]_{\im} = k[S']_{\im}$. Hence $k[S] = k[S']$.
\qed 

\brm
From the proof in this paper it is clear that the main result does not hold true for higher dimensional semigroup ring or even for semigroup rings of codimension two if the dimension of the ambient space is more than $4$. So the natural question is: Can one characterize all codimension two semigroup rings whose canonical module is Cohen-Macaulay. 
\erm

\end{document}